\documentclass[aos,MSNbibl,nameyear,dvips]{arximspdf}
\usepackage{graphicx}

\doi{10.1214/12-AOS1008} %kopijuoti is PTS
\volume{40}
\issue{2}
\pubyear{2012}
\firstpage{1263}
\lastpage{1282}

\makeatletter
\newtheorem{theorem}{Theorem}[section]
\newtheorem{corollary}[theorem]{Corollary}
\newcommand{\cov}{\operatorname{cov}}
\newcommand{\cor}{\operatorname{cor}}
\newcommand{\var}{\operatorname{var}}
\newcommand{\E}{\mathbb{E}}
\renewcommand{\citep}[1]{[\citeauthor{#1} (\citeyear{#1})]}
\newcommand{\citecc}[1]{\citeauthor{#1} (\citeyear{#1})}
\newcommand{\citeccs}[1]{\citeauthor{#1} [(\citeyear{#1})}
\newcommand{\eqref}[1]{(\ref{#1})}
\makeatother

\begin{document}
\begin{frontmatter}

\title{Rerandomization to improve covariate balance in experiments\thanksref{T1}}

\runtitle{Rerandomization}

\thankstext{T1}{Supported in part by the following Grants: NSF SES-0550887, NSF IIS-1017967 and NIH R01DA23879.}

\begin{aug}
\author[A]{\fnms{Kari Lock} \snm{Morgan}\corref{}\ead[label=e1]{kari@stat.duke.edu}}
\and
\author[B]{\fnms{Donald B.} \snm{Rubin} \ead[label=e2]{rubin@stat.harvard.edu}}

\runauthor{K. L. Morgan and D. B. Rubin}
\affiliation{Duke University and Harvard University}
\address[A]{Department of Statistics\\
Duke University\\
Box 90251\\
Durham, North Carolina 27708\\
USA\\
\printead{e1}} %adresu isvedimo komanda gale!
\address[B]{Statistics Department\\
Harvard University\\
1 Oxford St.\\
Cambridge, Massachusetts 02138\\
USA\\
\printead{e2}}
\end{aug}

% HISTORY:
\received{\smonth{11} \syear{2011}}
\revised{\smonth{4} \syear{2012}}

% ABSTRACT
\begin{abstract}
Randomized experiments are the ``gold standard'' for estimating causal
effects, yet often in practice, chance imbalances exist in covariate
distributions between treatment groups.  If covariate data are
available before units are exposed to treatments, these chance
imbalances can be mitigated by first checking covariate balance {\em
before} the physical experiment takes place.  Provided a precise
definition of imbalance has been specified in advance, unbalanced
randomizations can be discarded, followed by a rerandomization, and
this process can continue until a randomization yielding balance
according to the definition is achieved.  By improving covariate
balance, rerandomization provides more precise and trustworthy
estimates of treatment effects.\looseness=1
\end{abstract}

% KEYWORDS

\end{frontmatter}

%s1 ###
\section{A brief history of rerandomization} \label{intro}

Randomized experiments are the ``gold standard'' for estimating causal
effects, because randomization balances all potential confounding
factors {\em on average}.  However, if in a particular experiment, a
randomization creates groups that are notably unbalanced on important
covariates, should we proceed with the experiment, rather than
rerandomizing and conducting the experiment on balanced groups?

With $k$ independent covariates, the chance of {\em at least one}
covariate showing a ``significant difference'' between treatment and
control groups, at significance level $\alpha$, is $1 - (1-\alpha)^k$.
For a modest 10 covariates and a 5\% significance level, this
probability is 40\%.  ``Most experimenters on carrying out a random
assignment of plots will be shocked to find how far from equally the
plots distribute themselves'' \citep{fisher26}.  The danger of relying
on pure randomization to balance covariates has been described in
\citet{seidenfeld81}; \citet{urbach85}; \citet{krause03};
\citet{rosenberger08}; \citet{rubin08a}; \citet{keele09} and
\citet{worrall10}. Also, there exists much discussion historically over
whether randomization should be preferred over a purposefully balanced
assignment [\citecc{gossett38}; \citecc{yates39}; \citecc{greenberg51};
\citecc{harville75}; \citecc{arnold86}; \citecc{kempthorne86}]. Our
view is that with rerandomization, we can retain the advantages of
randomization, while also ensuring balance.

It is standard in randomized experiments today to collect covariate
data and check for covariate balance, yet typically this is done after
the experiment has started.  If covariate data are available before the
physical experiment has started, a randomization should be checked for
balance {\em before} the physical experiment is conducted.  If lack of
balance is noted, as Gosset stated, ``it would be pedantic to continue
with an arrangement of plots known beforehand to be likely to lead to a
misleading conclusion'' \citep{gossett38}.   It appears that Fisher
would agree.  In \citet{rubin08a}, Rubin recounts the following
conversation with his advisor Bill Cochran:

{\em Rubin}: What if, in a randomized experiment, the chosen randomized
allocation exhibited substantial imbalance on a prognostically
important baseline covariate?

{\em Cochran}: Why didn't you block on that variable?

{\em Rubin}: Well, there were many baseline covariates, and the correct
blocking wasn't obvious; and I was lazy at that time.

{\em Cochran}: This is a question that I once asked Fisher, and his
reply was unequivocal:

{\em Fisher (recreated via Cochran)}: Of course, if the experiment had
not been started, I would rerandomize.

\noindent A similar conversation between Fisher and Savage, wherein
Fisher advocates rerandomization when faced with an undesirable
randomization, is documented in \citeccs{savage62}, page 88].

Checking covariates and rerandomizing when needed for balance has been
advocated repeatedly.    \citet{sprott93} recommend
rerandomization when ``obvious'' lack of balance is observed.
\citet{rubin08a} suggests that if  ``important imbalances exist,
rerandomize, and continue to do so until satisfied.''  For clinical
trials,   \citet{worrall10} states that ``if such baseline
imbalances are found then the recommendation\,\ldots\ is to
re-randomize in the hope that this time no baseline imbalances will
occur.''   \citet{cox09} and   \citet{bruhn09} have advocated
rerandomization, suggesting either to do multiple randomizations and
pick the ``best,'' or to specify a~bound for the difference in
treatment and control covariate means for each covariate, following the
``Big Stick'' method of \citet{soares85}, and rerandomize until all
differences are within these bounds.  The latter rerandomization method
was used in \citet{maclure06}.

There are also many sources giving reasons {\em not} to rerandomize.
Good accounts of the debate over rerandomization can be found in
\citet{urbach85} and \citet{raynor86}. The most common critique of
rerandomization is that forms of analysis utilizing Gaussian
distribution theory are no longer valid [\citecc{fisher26};
\citecc{anscombe48}; \citecc{grundy50}; Holschuh\break (\citeyear{holschuh80});
\citecc{bailey83}; \citecc{urbach85}; \citecc{bailey86};
Bailey and Rowley (\citeyear{bailey87})]. Rerandomization changes the distribution of the
test statistic, most notably by decreasing the true standard error,
thus traditional methods of analysis that do not take this into account
will result in overly ``conservative'' inferences in the sense that
tests will reject true null hypotheses less often than the nominal
level and confidence intervals will cover the true value more often
than the nominal level. However, randomization-based inference is still
valid [\citecc{anscombe48}; \citecc{kempthorne55};
\citecc{brillinger78}; \citecc{tukey93}; \citecc{rosenberger02};
\citecc{moulton04}], because the rerandomization can be accounted for
during analysis.

All other critiques of rerandomization, of which we are aware, deal
with ``ad-hoc'' rerandomization, that is, rejecting randomizations
without specifying a rejection criterion in advance.  We only advocate
rerandomization if the decision to rerandomize or not is based on a
pre-specified criterion.  By specifying an objective rerandomization
rule before randomizing, and then analyzing results using
randomization-based methods, we can, in most circumstances, finesse all
existing criticisms of rerandomizing.

Some may think that rerandomization is unnecessary with large sample
sizes, because as the sample size increases, the difference in
covariate means between groups gets smaller, essentially proportional
to the square root of the sample size.  However, at the same rate,
confidence intervals and significance tests are getting more sensitive
to small differences in outcome means, which can be driven by small
differences in covariate means.

Despite the ongoing discussion about rerandomization and the fact that
it~is widely used in practice [\citecc{holschuh80}; \citecc{urbach85};
Bailey and Rowley (\citeyear{bailey87}); \citecc{imai08}; \citecc{bruhn09}], little has been
published on the mathematical implications of rerandomization.
Remarkably, it appears that no source even makes explicit the
conditions under which rerandomization is valid.  Although a few
rerandomization methods have been proposed [\citecc{moulton04};
\citecc{maclure06}; Bruhn and McKenzie (\citeyear{bruhn09}); \citecc{cox09}], the implications
have not been theoretically explored, to the best of our knowledge. The
only published theoretical results accompanying a rerandomization
procedure appear to be those in \citet{cox82}, which proposed
rerandomization to lower the sampling variance of covariance-adjusted
estimates.  Here we aim to fill these lacuna by (a) making explicit the
sufficient conditions under which rerandomization is valid, (b)
describing in detail a principled procedure for implementing
rerandomization and (c) providing corresponding theoretical results.

%f1 ###
\begin{figure}

\includegraphics{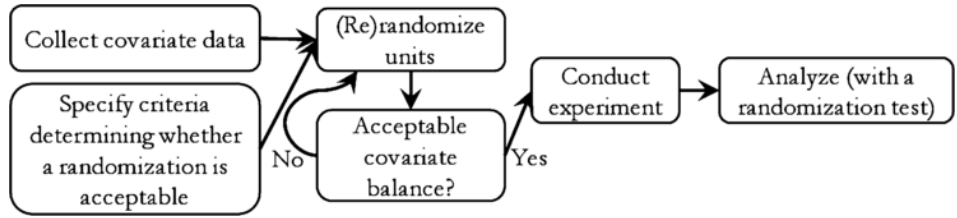}

 \caption{Procedure for implementing rerandomization.}
\label{flowchart}
\end{figure}

%s2 ###
\section{Rerandomization in general} \label{implementing}

%s2.1 ###
\subsection{Procedure} \label{procedure}

The procedure for implementing rerandomization is depicted in
Figure~\ref{flowchart}, and has the following steps:
\begin{longlist}
\item[(1)] Collect covariate data.
\item[(2)] Specify a balance
criterion determining when a randomization is acceptable.
\item[(3)]
Randomize the units to treatment groups.
\item[(4)] Check the balance
criterion; if the criterion is met, go to Step (5).  Otherwise, return
to Step (3).
\item[(5)] Conduct the experiment using the final
randomization obtained in Step (4).
\item[(6)] Analyze the results
using a randomization test, keeping only simulated randomizations that
satisfy the balance criterion specified in Step~(2).
\end{longlist}

Let $\mathbf{x}$ be the $n \times k$ covariate matrix representing $k$
covariates measured on $n$ experimental units.  Here we assume that a
sample of units has already been selected and is fixed.  Because we are
not considering the sampling mechanism, we are only interested in the
extent to which a causal effect estimate obtained in {\em this}
randomized experiment is a good estimate of the true causal effect
within the selected sample.  The $\mathbf{x}$ matrix includes all the
observed covariates for which balance between groups is desired, which
may include original covariates, and any functions of original
covariates, such as transformations, powers and interactions.  Let
$\mathbf{W}$ be the $n$-dimensional treatment assignment vector indicating
the treatment group for each unit.   The rerandomization criterion is
based on a row-exchangeable scalar function of~$\mathbf{x}$ and $\mathbf{W}$.
\[
\varphi  ( \mathbf{x},\mathbf{W}  ) =
\cases{ 1, &\quad \mbox{if $\mathbf{W}$ is an acceptable randomization,} \cr
0, &\quad\mbox{if $\mathbf{W}$ is not an acceptable randomization.}}
\]
The function $\varphi$ can vary depending on the relative importance of
balancing different covariates, on the level of covariate balance
desired and on the computational power available, but it is specified
in advance.

More generally, we can define a set of acceptance criteria,
$\mathcal{S}=\{ \varphi_s \}$,  from which we choose at each step, $s$,
either deterministically or stochastically, where this choice can
depend on the step, so that, for example, we can become more lenient as
the steps increase without success.  In this more general situation,
$\varphi_s  ( \mathbf{x}, \mathbf{W}  )$ denotes the acceptance criterion for step
$s$.  Although our theoretical results in Sections \ref{implementing} and \ref{mahalanobis} hold for
this more general setup, in practice we expect that the common choice
will be one function for all steps, and to avoid notational clutter, we
present results with one criterion.\looseness=1

Once $\varphi$ has been specified, units are randomized to treatment
groups (Step 3).  In the simplest form of rerandomization, this can be
done with no restrictions; for example, randomly choose an assignment
vector $\mathbf{W}$ from all possible vectors, or equivalently from all
possible partitions of the units into groups.  In practice, the initial
randomization is typically done with some restriction to equalize
treatment group sizes.

Rerandomization is simply a tool that allows us to draw from some
predetermined set of acceptable randomizations, $\{ \mathbf{W} \mid
\varphi (\mathbf{x}, \mathbf{W}) = 1\}$.    Rerandomization is analogous to
rejection sampling; a way to draw from a set that may be tedious to
enumerate.  Specifying a set of acceptable randomizations and then
choosing randomly from this set is recommended by
\citeauthor{kempthorne55} (\citeyear{kempthorne55,kempthorne86}) and \citet{tukey93},
and \citet{moulton04} notes that rerandomization may be required for
implementation of this idea when the set of acceptable randomizations
is difficult to enumerate a priori.

Within this framework, rerandomization simply generalizes classical
experimental designs.  For the basic completely randomized experiment
with fixed sample sizes in each treatment group, $\varphi (\mathbf{x}, \mathbf{W})
= 1$ when the number of units assigned to each group matches the
predetermined group sizes.  For a randomized block experiment, $\varphi
(\mathbf{x}, \mathbf{W}) = 1$ when predetermined numbers of units within each block
are assigned to each treatment group.  For a Latin square, $\varphi
(\mathbf{x}, \mathbf{W}) = 1$ when the randomization satisfies the Latin square
design.  These classical designs can be readily sampled from, so
rerandomization is computationally inefficient, although equivalent,
but for other functions,~$\varphi$, rerandomization may be a more
straightforward technique.  Rerandomization can also be used together
with any classical design.  For example, in a medical experiment on
hypertensive drugs, we may block on sex and a~coarse categorization of
baseline blood pressure, and use rerandomization to balance the
remaining covariates, including fine baseline blood
pressure.\looseness=1

Researchers are free to chose any function $\varphi$, provided it is
chosen in advance. Section \ref{unbiasedsec} describes the conditions
necessary to maintain general unbiasedness of simple point estimation,
Section \ref{mahalanobis} recommends a particular class of functions
and studies theoretical properties of this choice and Section~\ref{affineRR} discusses some reasons for choosing an affinely
invariant $\varphi$.

%s2.2 ###
\subsection{Analysis by randomization tests}\label{analysis}

Under most forms of rerandomization, increasing balance in the
covariates will typically create more precise estimated treatment
effects,\vadjust{\goodbreak} making traditional Gaussian distribution-based forms of
analysis statistically too conservative.  However, the final data can
be analyzed using a randomization test, maintaining valid frequentist
properties.  As Fisher stated, ``It seems to have escaped recognition
that the physical act of randomization\,\ldots\ affords the means, in
respect of any particular body of data, of examining a wider hypothesis
in which no normality of distribution is implied'' \citep{fisher35}.
This physical act of randomization need not be pure randomization, but
any randomization scheme that can be replicated when conducting the
randomization test.

We are interested in the effect of treatment assignment, $ \mathbf{W}$, on
an outcome, $\mathbf{y}$.  Let $y_i(W_i)$ denote the $i$th unit's, $\{ i =
1, \dots, n\}$, potential outcome under treatment assignment $W_i$,
following the Rubin causal model \citep{rubin74}.  Although
rerandomization can be applied to any number of treatment conditions,
to convey essential ideas most directly, we consider only two, and
refer to these conditions as treatment and control.  Let
\[
W_i =
\cases{ 1, &\quad \mbox{if treated,}
\cr
0, &\quad\mbox{if control.}
}
\]
Let $\mathbf{Y}_{obs( \mathbf{W})}$ denote the vector of observed outcome
values:
%e1 ###
\begin{equation}
Y_{obs, i} = y_i(1)W_i  +  y_i(0) (1-W_i),
\end{equation}
where for notational simplicity the subscript $obs$ means $obs( \mathbf{W})$.  Under the sharp null hypothesis of no treatment effect on any
unit, $y_i(1) = y_i(0)$ for every $i$, and thus the vector $\mathbf{Y}_{obs}$ is the same for every treatment assignment $ \mathbf{W}$.
Consequently, leaving $\mathbf{Y}_{obs}$ fixed and simulating many
acceptable randomization assignments, $\mathbf{W}$, we can empirically
create the distribution of any estimator, $g(\mathbf{x}, \mathbf{W}, \mathbf{y}_{obs})$, if
the null hypothesis were true.  To account for the rerandomization,
each simulated randomization must also satisfy $\varphi (\mathbf{x}, \mathbf{W}) =
1$.   Once the desired number of randomizations has been simulated, the
proportion of simulated randomizations with estimated treatment effect
as extreme or more extreme than that observed in the experiment is the
$p$-value.  Although a full permutation test (including all the
acceptable randomizations) is necessary for an exact $p$-value, the
number of simulated randomizations can be increased to provide a
$p$-value with any desired level of accuracy.  This test can incorporate
whatever rerandomization procedure was used, will preserve the
significance level of the test \citep{moulton04} and works for any
estimator. \citet{brillinger78}, \citet{tukey93} and
\citeccs{rosenberger02}, Chapter 7] suggest using randomization tests to
assess significance when restricted randomization schemes are used.

Because analysis by a randomization test requires generating many
acceptable randomizations, computational time can be important
to\vadjust{\goodbreak}
consider in advance.  Define $p_a \equiv P(\varphi = 1)$ to be the
proportion of acceptable randomizations.   The choice of $p_a$ involves
a trade-off between better balance and computational time; smaller
values of $p_a$ ensure better balance, but they also imply a longer
expected waiting time to obtain an acceptable randomization, at least
without clever computational devices.  The number of randomizations
required to get one acceptable randomization follows a geometric
distribution with parameter $p_a$, so $N$ simulated acceptable
randomizations for a randomization test will require on average $N/p_a$
randomizations to be generated.

The chosen $p_a$ must leave enough acceptable randomizations to perform
a randomization test.  In practice this is rarely an issue, because the
number of possible randomizations is huge even for modest $n$.  To
illustrate, the number of possible randomizations for $n = \{30, 50,
100\}$ randomizing to two equally sized treatment groups, ${n \choose
n/2}$, is on the order of $\{10^8, 10^{14}, 10^{29}\}$, respectively.
However, for small sample sizes, care should be taken to ensure the
number of acceptable randomizations does not become too small, for
example, less than 1000.

By employing the duality between confidence intervals and tests, for
additive treatment effects a confidence interval can be produced from a
randomization distribution as the set of all values for which the
observed data would not reject such a null hypothesized value
[\citecc{lehmann05}; \citecc{manly07}, Section~3.5, Section 1.4].
\citet{garthwaite96} provides an efficient algorithm for generating a
confidence interval for additive effects from a randomization test.
The assumption of additivity is statistically conservative, at least
asymptotically, as implied by Neyman's \citep{neyman23} results on
standard errors being overestimated when assuming it relative to the
actual standard errors.  A randomization test can be applied to any
sharp null hypothesis, that is, a hypothesis such that the observed
data implies specific values for all missing potential outcomes.

When the covariates being balanced are correlated with the outcome
variable, then rerandomization increases precision (Section
\ref{precisionsec}).  A randomization test reflects this increase in
precision.  Standard asymptotic-based frequentist analysis procedures
that do not take the rerandomization into account will be statistically
conservative.   Not only will distribution-based standard errors not
incorporate the increase in precision, but the act of rerandomizing
itself will increase the {\em estimated} standard error beyond that of
pure randomization.  If the total variance in the outcome is fixed,
decreasing the actual sampling variance between treatment group means
(by ensuring better balance), increases the variance within groups, and
it is this variance within groups that is traditionally used to
estimate the standard error \citep{fisher26}.  Thus, although
rerandomization decreases the {\em true} standard error, it actually
increases the standard error as estimated by traditional methods.  For
both of these reasons, the regular estimated standard errors will
overestimate the true standard error, and using the corresponding
distribution-based methods of analysis after rerandomization results in
overly wide confidence intervals and less powerful tests of hypotheses.

%s2.3 ###
\subsection{Maintaining an unbiased estimate} \label{unbiasedsec}

Although not needed to motivate rerandomization, we assume one goal is
to estimate the average treatment effect
\begin{eqnarray}\label{pace}
\tau &\equiv& \overline{y(1)} - \overline{y(0)}  \nonumber \\[-8pt]\\[-8pt]
&=& \frac{\sum_{i=1}^n y_i(1)}{n} - \frac{\sum_{i=1}^n y_i(0)}{n}.\nonumber
\end{eqnarray}
The fundamental problem in causal inference is that, because we only
observe $y_i(W_i)$ for each unit, we cannot calculate \eqref{pace}
directly, and we must estimate $\tau$ using only $\mathbf{Y}_{obs}$.  In
this section, we assume the Stable Unit Treatment Value Assumption
(SUTVA) \citep{rubin80}: the potential outcomes are fixed and do not
change with different possible assignment vectors $\mathbf{W}$.

The average treatment effect $\tau$ is usually estimated by the
difference in observed sample means,
\[\label{tauhat}
\hat{\tau}  \equiv \overline{Y}_{obs,T} - \overline{Y}_{obs,C},
\]
where
\[
\overline{Y}_{obs,T} \equiv \frac{\sum_{i=1} W_i y_i(1)}{\sum_{i=1}^n
W_i} \quad\mbox{and}\quad \overline{Y}_{obs, C} \equiv \frac{\sum_{i=1}^n
(1-W_i) y_i(0)}{\sum_{i=1}^n (1-W_i)}.
\]

\begin{theorem}\label{unbiased}
Suppose $\sum_{i=1}^n W_i = \sum_{i=1}^n (1-W_i)$ and
$\varphi(\mathbf{x},\mathbf{W}) = \varphi(\mathbf{ x},\allowbreak
\mathbf{1-W})$; then $\E(\hat{\tau} \mid \mathbf{ x}, \varphi =
1) = \tau$.
\end{theorem}

\begin{pf}
Under the specified conditions, $\mathbf{ W}$ and $\mathbf{ 1 - W}$ are
exchangeable.  Therefore, after rerandomization $\E  (W_i \mid \mathbf{
x},\varphi = 1  ) = \E  ( 1 - W_i \mid   \mathbf{ x},\varphi = 1  )$ $\forall
i$, so $\E  ( W_i \mid  \mathbf{ x},\varphi = 1  ) = \E  ( 1 - W_i \mid
\mathbf{ x},\varphi = 1  ) = 1/2$ $\forall i$. Hence
\begin{eqnarray*}
\E  (\hat{\tau} \mid \mathbf{ x}, \varphi = 1  )
&=& \E  \biggl( \frac{\sum_{i=1}^n W_i Y_{i, obs}}{n/2} - \frac{\sum_{i=1}^n (1-W_i) Y_{i, obs}}{n/2}  \Bigm|     \mathbf{ x},\varphi = 1  \biggr) \\
&=& \E  \biggl( \frac{\sum_{i=1}^n W_i y_i(1)}{n/2} - \frac{\sum_{i=1}^n (1-W_i)y_i(0)}{n/2}  \Bigm|    \mathbf{ x},\varphi = 1 \biggr) \\
&=& \frac{\sum_{i=1}^n \E(W_i \mid  \varphi = 1) y_i(1)}{n/2} - \frac{\sum_{i=1}^n  (1-\E (W_i \mid  \mathbf{ x},\varphi = 1 )  ) y_i(0)}{n/2}  \\
&=&\frac{\sum_{i=1}^n (1/2) y_i(1)}{n/2} - \frac{\sum_{i=1}^n (1/2) y_i(0)}{n/2}  \\
&=& \tau.
\end{eqnarray*}
\upqed\end{pf}

Theorem \ref{unbiased} holds for all outcome variables.  Corollary
\ref{covunbiased} follows by the same logic.

\begin{corollary}\label{covunbiased}
If $\sum_{i=1}^n W_i = \sum_{i=1}^n (1-W_i)$ and $\varphi(\mathbf{ x},
\mathbf{W}) = \varphi(\mathbf{ x},\break \mathbf{1-W})$, then
$\E(\overline{V}_T - \overline{V}_C  \mid \mathbf{ x}, \varphi = 1) =
0$ for any observed or unobserved covariate $V$.
\end{corollary}

If sample sizes are not fixed in advance, but each unit has $E(W_i \mid
\mathbf{ x}) = 1/2$ in the initial randomization, $\hat{\tau}$ is only
necessarily an unbiased estimate under the assumption of additivity. As
a small example under nonadditivity, consider $\mathbf{ x} = (0,1,2)$,
$\mathbf{ y}(1) = (1,1,0)$ and $\mathbf{ y}(0) = (0,0,1)$.  When
$\varphi(\mathbf{ x}, \mathbf{W}) = 1$ if the difference in $x$ means
between the two groups is $0$ and $\varphi = 0$ otherwise, the only two
acceptable randomizations are $\mathbf{ W} = (0,1,0)$ and $\mathbf{ W}
= (1,0,1)$.  For either acceptable randomization, $\hat{\tau} = 1/2$,
yet $\tau = 1/3$.  This artificial example also illustrates that if the
treatment groups are of unequal size, $\hat{\tau}$ will not necessarily
be an unbiased estimate after rerandomization.  If the treatment group
includes two units and the control group one unit, and
$\varphi(\mathbf{ x}, \mathbf{W})$ is the same as before, then the only
acceptable randomization is $\mathbf{ W} = (1,0,1)$, and once again,
$\hat{\tau} = 1/2$, whereas $\tau = 1/3$.

%s3 ###
\section{Rerandomization using Mahalanobis distance}\label{mahalanobis}

To simplify the statement of theoretical results, we assume the sample
sizes for the treatment and control groups are fixed in advance, with
$p_w$ the fixed proportion of treated units,
%e2 ###
\begin{equation}
p_w = \frac{\sum_{i=1}^n W_i}{n}.
\end{equation}

Let $\overline{\mathbf{ X}}_T - \overline{\mathbf{ X}}_C$ be the
$k$-dimensional vector of the difference in covariate means between the
treatment and control groups,
%e3 ###
\begin{equation}\label{diffxmeansW}
\overline{\mathbf{ X}}_T - \overline{\mathbf{ X}}_C  = \frac{\mathbf{
x' W}}{np_w} - \frac{\mathbf{ x'(1-W)}}{n(1-p_w)} = \frac{\mathbf{
x}'(\mathbf{ W} - p_w\mathbf{ 1})}{np_w(1-p_w)}.
\end{equation}

We consider Mahalanobis distance as a rerandomization criterion because
it is an affinely invariant scalar measure of multivariate covariate
balance.  Mahalanobis distance is defined by
%e5 ###
%e4 ###
\begin{eqnarray}
\label{M1}M &\equiv&  (\overline{\mathbf{ X}}_T - \overline{\mathbf{
X}}_C  )'  [\cov  (\overline{\mathbf{ X}}_T - \overline{\mathbf{ X}}_C
)  ]^{-1}
(\overline{\mathbf{ X}}_T - \overline{\mathbf{ X}}_C  )  \\
\label{M}&=& np_w(1-p_w) ( \overline{\mathbf{ X}}_T -
\overline{\mathbf{ X}}_C  )' \cov(\mathbf{x})^{-1}  (
\overline{\mathbf{ X}}_T - \overline{\mathbf{ X}}_C
 ),
\end{eqnarray}
where $\cov(\mathbf{ x})$ represents the sample covariance matrix of $\mathbf{
x}$.  The quantities~$n$, $p_w$ and $\cov(\mathbf{ x})$ are known and
constant across randomizations.  If $\cov(\mathbf{ x})$ is singular, for
example, if $k \geq n$, then $\cov(\mathbf{ x})^{-1}$ can be replaced with
the pseudo-inverse of $\cov(\mathbf{ x})$.  For cluster randomized
experiments, see \citet{hansen08}.

Due to the finite population central limit theorem, $\overline{\mathbf{
X}}_T - \overline{\mathbf{ X}}_C$ is asymptotically multivariate
normally distributed over its randomization distribution
[\citecc{erdos59}; \citecc{hajek60}]. Normality of $\overline{\mathbf{
X}}_T - \overline{\mathbf{ X}}_C$ is not necessary for rerandomization,
but assuming normality allows for the theoretical results of this
section. If $\overline{\mathbf{ X}}_T - \overline{\mathbf{ X}}_C$ is
multivariate normal, then under pure randomization, $M \sim \chi_k^2$
[\citet{mardia80}, page 62]; $M$ is the statistic used in Hotelling's
$T^2$ test, but note that here $M$ follows a $\chi_k^2$ distribution
because $\mathbf{ x}$ is considered fixed.

A randomization is deemed ``acceptable'' whenever $M$ falls below a
certain threshold, $a$. Let $p_a$ be the proportion of randomizations
that are acceptable, so that $P(M \leq a) = p_a$.  Either $a$ or $p_a$
can be specified a priori, and then the other is fixed either using $M
\sim \chi_k^2$ if sample sizes are large enough or using an empirical
distribution of $M$ achieved through simulation.   The rerandomization
criterion, $\varphi_M$, is
%e6 ###
\begin{equation}
\varphi_M  ( \mathbf{ x},\mathbf{ W}  ) \equiv \cases{ 1, &\quad
\mbox{if $M \leq a$,} \cr 0 ,&\quad\mbox{if $M > a$.} }
\end{equation}

%s3.1 ###
\subsection{Covariate balance under $\varphi_M$}\label{balance}

\begin{theorem} \label{covariates}
Assume rerandomization is conducted using $\varphi_M$ with $p_w = 1/2$,
and the covariate means are multivariate normal; then
%e7 ###
\begin{equation}\label{covX}
\cov  (\overline{\mathbf{ X}}_T - \overline{\mathbf{ X}}_C  \mid
\mathbf{ x}, \varphi_M = 1  ) = v_a  \cov  (\overline{\mathbf{ X}}_T -
\overline{\mathbf{ X}}_C \mid \mathbf{ x}, ),
\end{equation}
where
%e8 ###
\begin{equation}\label{va}
v_a \equiv \frac{2}{k} \times \frac{\gamma  (k/2 +1, a/2
)}{\gamma  (k/2, a/2  )} = \frac{P ( \chi_{k+2}^2 \leq a
)}{P ( \chi_k^2 \leq a  )}
\end{equation}
and $\gamma$ denotes the incomplete gamma function: $\gamma(b,c) \equiv
\int_0^c y^{b-1} e^{-y} \,dy$.
\end{theorem}

The proof of Theorem \ref{covariates} is in the
\hyperref[app]{Appendix}.

\begin{figure}

\includegraphics{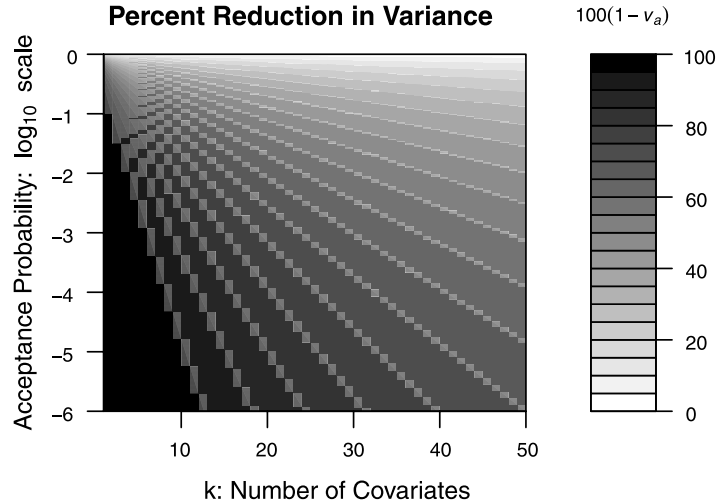}

\caption{The percent reduction in variance for each
covariate difference in means, as a function of the number of
covariates and the proportion of randomizations accepted.}
\label{privfig}
\end{figure}

In the field of matching, emphasis has been placed on ``percent
reduction in bias'' \citep{cochran73}.  In the context of randomized
experiments there is no bias, and rerandomization instead reduces the
sampling variance of the difference in covariate means, yielding
differences that are more closely concentrated around $0$. Define the {\em percent reduction in variance}, the percentage by which
rerandomization reduces the randomization variance of the difference in
means, for each covariate, $x_j$, by
%e9 ###
\begin{equation}\label{priv}
100  \biggl(\frac {\var  ( \overline{X}_{j, T} - \overline{X}_{j, C}
\mid \mathbf{ x}  ) - \var  ( \overline{X}_{j, T} - \overline{X}_{j, C}
\mid \mathbf{ x}, \varphi = 1 )}  {\var  ( \overline{X}_{j, T} -
\overline{X}_{j, C}  \mid \mathbf{ x}  )}  \biggr).
\end{equation}
By Theorem \ref{covariates}, the percent reduction in variance for each
covariate, and for any linear combination of these covariates, is
%e10 ###
\begin{equation}\label{privM}
100(1-v_a)
\end{equation}
and is shown as a function of $k$ and $p_a$ in Figure~\ref{privfig},
where by \eqref{va}, $0 \leq v_a \leq 1$.  The lower the acceptance
probability and the fewer covariates being balanced, the larger the
percent reduction in variance.

%f2 ###

Notice that Theorem \ref{covariates} holds for any covariate
distribution, as long as the sample size is large enough for the
central limit theorem to ensure normally distributed covariate means.
An exact value is not needed, and an estimate is used only to guide the
choice of $p_a$; it has no influence on the validity of resulting
inferences.

%s3.2 ###
\subsection{Precision of the estimated treatment effect} \label{precisionsec}

Rerandomization improves precision, provided the outcome and covariates
are correlated.  Thus researchers can increase the power of tests and
decrease the width of confidence intervals simply at the expense of
computational time.

\begin{theorem} \label{precisionThm}
If \textup{(a)} rerandomization is conducted using $\varphi_M$ with $p_w = 1/2$,
\textup{(b)} the covariate and outcome means are normally distributed, and \textup{(c)}
the treatment effect is additive, then the percent reduction in
variance of $\hat{\tau}$~is
%e11 ###
\begin{equation}\label{privtau}
100(1-v_a) R^2,
\end{equation}
where $R^2$ represents the squared multiple correlation between
$\mathbf{y}$ and $\mathbf{x}$ within a treatment group and $v_a$ is as defined
in \eqref{va}.
\end{theorem}

\begin{pf}
Regardless of the true relationship between the outcome and covariates,
by additivity we can write
%e12 ###
\begin{equation} \label{truth}
y_i(W_i) \mid  \mathbf{ x}_i = \beta_0 + \beta '\mathbf{ x}_i + \tau W_i + e_i,
\end{equation}
where $\beta_0 + \beta'\mathbf{ x}_i$ is the projection of $y_i$ onto the
space spanned by $(\mathbf{ 1, x})$, and $e_i$ is a residual that
encompasses any deviations from the linear model.  Then the estimated
treatment effect, $\hat{\tau}$, can be expressed as
%e13 ###
\begin{equation}
\hat{\tau} = \beta'  ( \overline{\mathbf{ X}}_T - \overline{\mathbf{ X}}_C
 ) + \tau +  ( \overline{e}_T - \overline{e}_C  ).
\end{equation}
Because $\tau$ is constant and the first and last terms are
uncorrelated, we can express the variance of $\hat{\tau}$ as
\begin{eqnarray}
\var({\hat{\tau}}) &=&
\var  \bigl( \beta'  ( \overline{\mathbf{ X}}_T - \overline{\mathbf{ X}}_C  )  \bigr) + \var  (\overline{e}_T - \overline{e}_C  ) \nonumber\\[-8pt]\\[-8pt]
&=& \beta' \cov  ( \overline{\mathbf{ X}}_T - \overline{\mathbf{ X}}_C  ) \beta
+ \var  (\overline{e}_T - \overline{e}_C  ).\nonumber
\end{eqnarray}
By Theorem \ref{covariates}, rerandomization modifies the first term
  by the factor $v_a$.  Because under normality,
orthogonality implies independence, the difference in residual means is
independent of the difference in covariate means, and thus
rerandomization has no affect on the second term.  Therefore, the
variance of $\hat{\tau}$ after rerandomization restricting $M \leq a$
is
\begin{eqnarray}\label{vartauhat1}
\quad\var({\hat{\tau}}\mid \mathbf{ x}, M \leq a) &=& \beta' \cov  ( \overline{\mathbf{ X}}_T -
\overline{\mathbf{ X}}_C\mid \mathbf{ x},  M \leq a  ) \beta + \var  (\overline{e}_T - \overline{e}_C \mid \mathbf{ x},  M \leq a  ) \nonumber \\[-4pt]\\[-12pt]
&=&v_a  \beta' \cov  ( \overline{\mathbf{ X}}_T - \overline{\mathbf{ X}}_C \mid
\mathbf{ x}  ) \beta + \var  (\overline{e}_T - \overline{e}_C \mid \mathbf{ x}
 ).\nonumber
\end{eqnarray}
Let $\sigma_e^2$ be the variance of the residuals and $\sigma_y^2$ be
the variance of the outcome within each treatment group, where
$\sigma_e^2 = \sigma_y^2 (1-R^2)$.  Thus
%e14 ###
\begin{equation}\label{evar}
 \var  (\overline{e}_T - \overline{e}_C \mid \mathbf{ x}  ) = \frac{\sigma_e^2}{np_w(1-p_w)} = \frac{ \sigma_y^2 (1-R^2)}{np_w(1-p_w)},
\end{equation}
and
%e15 ###
\begin{eqnarray}\label{xvar}
\beta' \cov  ( \overline{\mathbf{X}}_T - \overline{\mathbf{X}}_C \mid \mathbf{ x}  )\beta
&=& \var(\hat{\tau} \mid \mathbf{x}) - \var  (\overline{e}_T - \overline{e}_C \mid \mathbf{ x}  ) \nonumber \\
&=&\frac{\sigma_y^2 - \sigma_e^2}{np_w(1-p_w)}   \\
&=& \frac{\sigma_y^2 R^2}{np_w(1-p_w)}.\nonumber
\end{eqnarray}

%f3 ###
\begin{figure}

\includegraphics{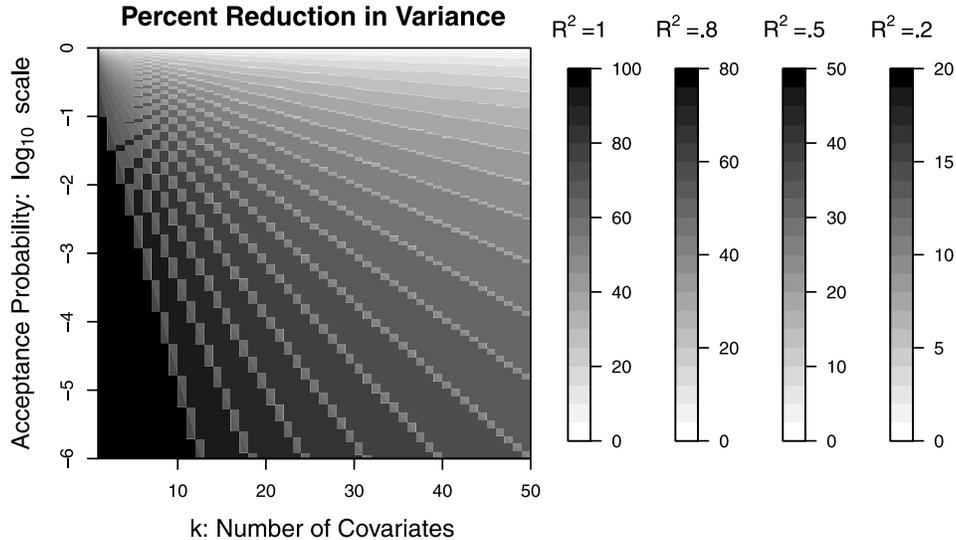}

\caption{The percent
reduction in variance for the estimated treatment effect, as a function
of the acceptance probability, the number of covariates, and $R^2$.}
\label{privy}
\end{figure}

Therefore by \eqref{vartauhat1}, \eqref{evar} and \eqref{xvar}, the
variance of $\hat{\tau}$ after rerandomization~is
\begin{eqnarray*}
\var(\hat{\tau} \mid  \mathbf{x}, M\leq a)
&=& v_a  \beta' \cov  ( \overline{\mathbf{X}}_T - \overline{\mathbf{X}}_C \mid \mathbf{x}  ) \beta + \var  (\overline{e}_T -
\overline{e}_C \mid \mathbf{x}  )   \\
&=& \frac{v_a\sigma_Y^2R^2}{np_w(1-p_w)} + \frac{\sigma_Y^2(1-R^2) }{np_w(1-p_w)}  \\
&=&   \bigl(1-(1-v_a)R^2 \bigr) \frac{\sigma_Y^2 }{np_w(1-p_w)}   \\
&= &  \bigl(1-(1-v_a)R^2 \bigr) \var(\hat{\tau} \mid \mathbf{x}).
\end{eqnarray*}
Thus the percent reduction in variance is $100( 1- (1-(1-v_a)R^2) =
100(1-v_a)R^2$.\vadjust{\goodbreak}
 \end{pf}

The percent reduction in variance for the estimated treatment effect,
shown as a function of $k$, $p_a$ and $R^2$ in Figure \ref{privy}, is
simply the percent reduction in variance for each covariate, scaled by
$R^2$.  Because under the specified conditions $\hat{\tau}$ is unbiased
by Theorem \ref{unbiased}, $100(1-v_a)R^2$ is not only the percent
reduction in variance in the estimated treatment effect, but also the
percent reduction in mean square error (MSE).

If regression (i.e., analysis of covariance) is used to adjust for
imbalance in a completely randomized experiment, the percent reduction
in variance~is
%e16 ###
\begin{equation}\label{ancova}
100  \biggl[ \biggl(1 + \frac{M}{n} \biggr)R^2 -
\frac{M}{n} \biggr]
\end{equation}
\citep{cox82}, where $M$ is as in \eqref{M}.  Comparing \eqref{ancova}
to \eqref{priv}, we see that rerandomization can increase precision
more than regression adjustment because there is no estimation of
regression coefficients with the former.  Note that the highest percent
reduction in variance achievable by either rerandomization or
regression is $100R^2$, achieved with perfect covariate mean balance.

%s4 ###
\section{Affinely invariant rerandomization criteria}\label{affineRR}

In this section we explore the theoretical implications of choosing an
affinely invariant rerandomization criterion, meaning that for any
affine transformation of $\mathbf{x}$, $a + \mathbf{b x}$,
$\varphi(\mathbf{x},\mathbf{W}) = \varphi(a + \mathbf{bx},
\mathbf{W})$. Measures based on inner products, such as Mahalanobis
distance or the estimated best linear discriminant, are affinely
invariant, as are criteria based on propensity scores estimated by
linear logistic regression \citep{rubin92}.  Results in this section
parallel those for affinely invariant matching methods \citep{rubin92}.

In the previous sections, we regarded $\mathbf{x}$ as fixed, and only
the randomization vector, $\mathbf{W}$ was random. In this section, to
use ellipsoidal symmetry of~$\mathbf{x}$, we regard both
$\mathbf{x}$ and $\mathbf{W}$ as random, so expectations are over
repeated draws of $\mathbf{x}$ and repeated randomizations.

\begin{theorem}\label{affineThm}
If $\varphi$ is affinely invariant, and if $\mathbf{x}$ is ellipsoidally
symmetric, then
%e18 ###
%e17 ###
\begin{eqnarray}
\E  ( \overline{\mathbf{X}}_T - \overline\mathbf{X}_C \mid \varphi = 1
 ) & =& \E  ( \overline{\mathbf{X}}_T - \overline\mathbf{X}_C  ) =
\mathbf{0} \quad \mbox{and}
\\
\cov  ( \overline{\mathbf{X}}_T - \overline\mathbf{X}_C \mid \varphi = 1
 ) &\propto & \cov  ( \overline{\mathbf{X}}_T - \overline\mathbf{X}_C
 ).
\end{eqnarray}
\end{theorem}

\begin{pf}
First, by ellipsoidal symmetry there is an affine transformation of
$\mathbf{x}$ to a canonical form with mean (center) zero and covariance
(inner product) $\mathbf{I}$, the $k$-dimensional identity matrix.  The
distribution of the matrix $\mathbf{x}$ in the treated group of size
$np_w$ and the control group of size $n(1-p_w)$ are both independent
and identically distributed samples from this zero centered spherical
distribution.  Any affinely invariant rule for selecting subsets of
treated and control units will be a function of affinely invariant
statistics in the treatment and control groups that are also
zero-centered spherically symmetric.  Applying $\varphi$ creates
concentric zero-centered sphere(s) that partition the space of these
statistics into regions where $\varphi = 1$ and $\varphi = 0$, and
therefore the distribution of such statistics remains zero-centered and
spherically symmetric.  Transforming back to the original form
completes the proof.~%
\end{pf}

\begin{corollary} \label{linearfunction}
If $\varphi$ is affinely invariant and if $\mathbf{x}$ is ellipsoidally
symmetric, then rerandomization leads to unbiased estimates of any
linear function of $\mathbf{x}$.
\end{corollary}

Note that, unlike Corollary \ref{covunbiased}, Corollary
\ref{linearfunction} applies no matter how the sample sizes are chosen.

\begin{corollary}\label{correlations}
If $\varphi$ is affinely invariant and if $\mathbf{x}$ is ellipsoidally
symmetric, then
%e19 ###
\begin{equation}
\cor  ( \overline{\mathbf{X}}_T - \overline\mathbf{X}_C \mid \varphi = 1
 )  =\cor  ( \overline{\mathbf{X}}_T - \overline\mathbf{X}_C  ).
\end{equation}
\end{corollary}

One possible method of rerandomization, suggested by \citet{moulton04},
\citet{maclure06}, \citet{bruhn09} and \citet{cox09}, is to place
bounds separately on each entry of $\overline{\mathbf{X}}_T -
\overline{\mathbf{X}}_C$ and ensure that each covariate difference is
within its specified caliper.  However, this method is not affinely
invariant and will generally destroy the correlational structure of
$\overline{\mathbf{X}}_T - \overline{\mathbf{X}}_C$, even when $\mathbf{x}$ is
ellipsoidally symmetric.

 Analogous to
``Equal Percent Bias Reducing'' (EPBR) matching methods
\citep{rubin76}, a rerandomization method is said to be ``Equal Percent
Variance Reducing'' (EPVR) if the percent reduction in variance is the
same for each covariate.

\begin{corollary}
If $\varphi$ is affinely invariant and if $\mathbf{x}$ is ellipsoidally
symmetric, then rerandomization is EPVR for $\mathbf{x}$ and any linear
function of $\mathbf{x}$.
\end{corollary}

Rerandomization methods that are not affinely invariant could increase
the variance of some linear combinations of covariates \citep{rubin76}.

Although affinely invariant methods have desirable properties in
general, they are not always preferred.  For example, if covariates are
known to vary in importance, a rerandomization method that is not EPVR
may be more desirable, allowing greater percent reduction in variance
for more important covariates.  Rerandomization criteria that take into
account covariates of varying importance are discussed in
\citeccs{lock11}, Chapter 4].

%s5 ###
\section{Discussion}\label{s5}

%s5.1 ###
\subsection{Alternatives for balancing covariates}

Rerandomization is certainly not the only way to balance covariates
before the experiment.

With only a few categorical covariates, simple blocking can
successfully balance all covariates, and there is no need for
rerandomization.  With many covariates each taking on many values,
however, blocking on all covariates can be impossible, and in this case
we recommend  blocking on the most important covariates, and
rerandomizing to balance the components of the covariates orthogonal to
the blocks.  Blocking and rerandomization can, and we feel should, be
used together.  Multivariate matching [\citecc{greevy04};
\citecc{rubin06}; \citecc{ho07}; \citecc{imai09}; Xu and Kalbfleisch (\citeyear{xu10})] is a
special case of blocking that can better handle many covariates.

Restricted (or constrained) randomization [\citecc{yates48};
Grundy and Healy (\citeyear{grundy50}); \citecc{youden72}; \citecc{bailey83}] restricts the
set of acceptable randomizations in a way that preserves the validity
of asymptotic-based distributional methods of analysis.  However, most
work on restricted randomization is specific to agricultural plots, and
apparently has not been extended to multiple covariates.  Blocking,
matching and restricted randomization can all also be implemented
through rerandomization by specifying the set of acceptable
randomizations through $\varphi$.

The Finite Selection Model (FSM) [\citecc{morris79}; \citecc{morris00}]
provides balance for multiple covariates, but provides a fixed amount
of balance in a\vadjust{\goodbreak} fixed amount of computational time.  Rerandomization
has the flexibility to choose the desired tradeoff between balance and
computational time.  More details comparing FSM with rerandomization
are in [\citet{lock11}, Section 5.5].

Covariate-adaptive randomization schemes [\citecc{efron71};
White and Freedman (\citeyear{white78}); \citecc{pocock75}; \citecc{pocock79};
\citecc{simon79}; \citecc{birkett85}; \citecc{aickin01};
\citecc{atkinson02}; \citecc{scott02}; \citecc{mcentegart03};
\citecc{rosenberger08}] are designed for clinical trials with
sequential treatment allocation over extended periods of time.
Rerandomization as proposed here is not applicable to sequential
allocation, and instead readers interested in such trials can refer to
the above sources.

If covariates are not balanced before the experiment, post-hoc methods
such as regression adjustment are commonly used, which rely on
assumptions that often cannot be verified [\citecc{tukey93};
\citecc{freedman08}]. Moreover, unlike post-hoc methods,
rerandomization is conducted entirely at the design stage, and so
cannot be influenced by outcome data. \citet{tukey93} and
\citet{rubin08b} give convincing reasons for why as much as possible
should be done in the design phase of an experiment, before outcome
data are available, rather than in the analysis stage when the
researcher has the potential to bias the results, consciously or
unconsciously.

%s5.2 ###
\subsection{Extensions and additional considerations}

For multiple treatment groups, any of the test statistics commonly used
in multivariate analysis of variance (MANOVA) can be used to measure
balance.  The standard statistics are all equivalent to Mahalanobis
distance in the special case of two groups.  Extensions for multiple
treatment groups are discussed in \citeccs{lock11}, Section 5.2].

For unbiased estimates using rerandomization with treatment groups of
unequal sizes, multiple treatment groups of equal size can be created,
and then merged as needed after the rerandomization procedure, but
before the physical experiment.  If extra units are discarded to form
equal sized treatment groups and rerandomization is employed, precision
can actually increase if covariates are highly correlated with the outcome
[\citet{lock11}, Section 5.3].

In a Bayesian analysis, as long as all covariates relevant to
$\varphi(\mathbf{x},\mathbf{W})$ are conditioned on, the design is ignorable
\citep{rubin78}, and theoretically, the analysis can proceed as usual.

%s6 ###
\section{Conclusion}

Randomization balances covariates across treatment\break groups, but only on
average, and in any one experiment covariates may be unbalanced.
Rerandomization provides a simple and intuitive way to improve
covariate balance in randomized experiments.

To perform rerandomization, a criterion determining whether a
randomization is acceptable needs to be specified.  For unbiasedness,
this rule needs to be\vadjust{\goodbreak} symmetric regarding the treatment groups.  If the
criterion is affinely invariant, then for ellipsoidally symmetric
distributions, balance improvement will be the same for all covariates
(and all linear combinations of the covariates), and correlations
between covariate differences in means will be maintained.  One such
criterion is to rerandomize whenever Mahalanobis distance exceeds a
certain threshold.

When the covariates are correlated with the outcome, rerandomization
increases precision.  If the analysis reflects the rerandomization
procedure, this leads to more precise estimates, more powerful tests
and narrower confidence intervals.

 %A randomization test, following the
%rerandomization procedure when generating each simulated randomization,
%will give $p$-values and confidence intervals with the nominal type I
%error and coverage rates.

\begin{appendix}\label{app}
\section*{Appendix}

\begin{pf*}{Proof of Theorem \ref{covariates}}
Because $M \sim \chi_k^2$ under pure randomization when the covariate
means are normally distributed, rerandomization affects the mean of $M$
as follows:
%e20 ###
\begin{eqnarray}
\E(M \mid \mathbf{x}, M \leq a) &=& \frac{ (1/(\Gamma(k/2)
2^{k/2})) \int_0^{a} y^{k/2} e^{-y/2}\, dy}
{(1/(\Gamma(k/2) 2^{k/2})) \int_0^{a} y^{k/2-1} e^{-y/2} \,dy} \nonumber \\
&=& \frac{\int_0^{a} (y/2  )^{k/2} e^{-y/2} \,dy}{(1/2) \int_0^{a} (y/2  )^{k/2-1} e^{-y/2} \,dy}   \\
&=& 2 \times \frac{\gamma  (k/2 +1, a/2  )}{\gamma
 (k/2, a/2  )}.\nonumber
\end{eqnarray}

To prove \eqref{covX}, we convert the covariates to canonical form
[Rubin and Thomas (\citeyear{rubin92})].  Let $\bolds{\Sigma} = \cov(\overline{\mathbf{X}}_T -
\overline{\mathbf{X}}_C \mid \mathbf{x})$, and define
%e21 ###
\begin{equation}\label{Z}
\mathbf{Z} \equiv \bolds{\Sigma}^{-1/2} (\overline{\mathbf{X}}_T -
\overline{\mathbf{X}}_C  ),
\end{equation}
where $\bolds{\Sigma}^{-1/2}$ is the Cholesky square root of
$\Sigma^{-1}$, so ${\bolds{\Sigma}^{-1/2}}'\bolds{\Sigma}^{-1/2}=
\bolds{\Sigma}^{-1}$. By the assumption of normality,
\[
\mathbf{Z} \mid \mathbf{x} \sim N_k (\mathbf{0}, \mathbf{I}  ).
\]
Due to normality, uncorrelated implies independent and thus the
elements of $\mathbf{Z}$ are independent and identically distributed
(i.i.d.) standard normals.  Therefore, the elements of $\mathbf{Z}$ are
exchangeable.

By \eqref{M1}, $M = \mathbf{Z}'\mathbf{Z} = \sum_{j=1}^k Z_j^2$.  Therefore
for each $j$ we have
%e23 ###
%e22 ###
\begin{eqnarray}
\var(Z_j \mid  \mathbf{x}, M \leq a) &=& \E(Z_j^2 \mid  \mathbf{x}, M \leq a) \nonumber \\
\label{exchangeable}&=& \frac{\E(M \mid  \mathbf{x}, M \leq a)}{k}  \\
&=& \frac{2}{k} \times \frac{\gamma  (k/2 +1, a/2  )}{\gamma  (k/2, a/2  )} \nonumber \\
&=& v_a,
\end{eqnarray}
where \eqref{exchangeable} follows from the exchangeability of the\vadjust{\goodbreak}
elements of $\mathbf{Z}$.

After enforcing $M\leq a$, the elements of $\mathbf{Z}$ are no longer
independent (they will be negatively correlated in magnitude), but with
signs they remain uncorrelated due to symmetry:
%e26 ###
%e25 ###
%e24 ###
\begin{eqnarray}
\cov(Z_i, Z_j \mid  \mathbf{x}, M\leq a)
&=& \E(Z_i Z_j \mid  \mathbf{x}, M\leq a) \nonumber
\\
&&{}- \E(Z_i \mid  \mathbf{x}, M\leq a)\E(Z_j \mid  \mathbf{x}, M \leq a) \nonumber \\
\label{unbiasedz}&=& \E \bigl( \E(Z_i Z_j \mid  Z_j, \mathbf{x}, M \leq a) \mid  \mathbf{x}, M \leq a  \bigr) - 0  \\
&=& \E \bigl( Z_j \E(Z_i \mid  Z_j, \mathbf{x}, M \leq a)\mid \mathbf{x}, M \leq a  \bigr)  \nonumber \\
\label{cov=0}&=&  \E ( Z_j \times 0\mid \mathbf{x}, M \leq a  )  \\
\label{covZ0}&=& 0,
\end{eqnarray}
where \eqref{unbiasedz} follows from Corollary \ref{covunbiased}, and
\eqref{cov=0} follows because $(Z_i \mid  Z_j, M \leq a) \sim (-Z_i
\mid  Z_j , M \leq a)$, thus $\E(Z_i \mid  Z_j, M \leq a)=0$ for all
$i, j$.

Thus after rerandomization the covariance matrix of $\mathbf{Z}$ is $v_a
\mathbf{I}$, hence
\begin{eqnarray*}
\cov(\overline{\mathbf{X}}_T - \overline{\mathbf{X}}_C\mid \mathbf{x}, M \leq a)
&=& \cov(\bolds{\Sigma}^{1/2}\mathbf{Z}\mid \mathbf{x}, M \leq a)    \\
&=& \bolds{\Sigma}^{1/2} \cov(\mathbf{Z} \mid \mathbf{x}, M \leq a)   {\bolds{\Sigma}^{1/2}}{}'   \\
&=& v_a \bolds{\Sigma}   \\
&=& v_a \cov (\overline{\mathbf{X}}_T - \overline{\mathbf{X}}_C \mid
\mathbf{x} ).
\end{eqnarray*}
\upqed\end{pf*}

\end{appendix}

\section*{Acknowledgments}

We appreciate the extraordinarily helpful comments of the editor,
Professor B\"uhlmann, and two reviewers.

% imsref loaded by arune.pranskunaite, 2012-05-31 08:48:46

\printaddresses

\end{document}